# Functional asymptotic confidence intervals for a common mean of independent random variables

**Yuliya V. Martsynyuk***

*School of Mathematics and Statistics, Carleton University,*
*1125 Colonel By Drive, Ottawa, ON K1S 5B6, Canada*
*e-mail:* ymartsynyuk@yahoo.com

**Abstract:** We consider independent random variables (r.v.'s) with a common mean $\mu$ that either satisfy Lindeberg's condition, or are symmetric around $\mu$. Present forms of existing functional central limit theorems (FCLT's) for Studentized partial sums of such r.v.'s on $D[0,1]$ are seen to be of some use for constructing asymptotic confidence intervals, or what we call functional asymptotic confidence intervals (FACI's), for $\mu$. In this paper we establish completely data-based versions of these FCLT's and thus extend their applicability in this regard. Two special examples of new FACI's for $\mu$ are presented.

## Contents

## 1. Introduction and main results

Estimating an unknown mean of a population has been a prominent classical problem in statistics. Perhaps, the most famous and influential work on this

*Research supported by NSERC Canada Discovery Grants of M. Csörgő and B. Szyszkowicz at Carleton University, and an NSERC Postdoctoral Fellowship of Yu. V. Martsynyuk at University of Ottawa.

subject is "Student" (1908) that celebrated its centennial last year. Facing the problem of an accurate interval estimation of a mean of a small random sample drawn from a normally distributed population with an unknown variance, W.S. Gosset ("Student"), among other things, concluded the exact distribution of the $1/\sqrt{n-1}$ multiple of what is now known as the Student statistic, or the so-called $t-$random variable with $n-1$ degrees of freedom.

Estimating a common mean of several populations has also been an outstanding statistical problem. It is frequently posed in the context of finite samples. A considerable amount of the literature in this regard treats the case of several normal populations with unknown and possibly unequal variances (cf., e.g., Graybill and Deal (1959), Normwood and Hinkelmann (1977), Pal and Kim (1997), and further references in these papers).

The present paper deals with asymptotic confidence interval estimation of a common mean of unspecified populations. Let $\{Z_i, i \geq 1\}$ be a sequence of independent, but not necessarily indentically distributed, random variables (r.v.'s) with a common mean $\mu$. We will consider two kinds of such r.v.'s, those with finite positive variances that satisfy Lindeberg's condition, and $Z_i$'s that are symmetric around $\mu$ and do not necessarily have finite variances.

The studies of the present paper were motivated in part by the problems the author faced in the context of linear error-in-variables models, when establishing functional asymptotic confidence intervals for the slope in such models in Martsynyuk (2008).

### *1.1. Review of invariance principles for Student processes based on independent random variables*

*Case of random variables satisfying Lindeberg's condition*

Suppose that $\mu = 0$ and $0 < \operatorname{Var} Z_i = \sigma_i^2 < \infty$, $i \geq 1$. Consider the Student statistic

$$T_n(Z_1, \ldots, Z_n) := \frac{\sum_{i=1}^n Z_i/\sqrt{n}}{(\sum_{i=1}^n (Z_i - \overline{Z})^2/(n-1))^{1/2}}, \tag{1}$$

with $\overline{Z} := n^{-1} \sum_{i=1}^n Z_i$, $n \geq 1$. In view of (1), one can define a Student process in $D[0,1]$ space as follows:

$$T_n^t(Z_1, \ldots, Z_n) := \frac{\sum_{i=1}^{K_n(t)} Z_i/\sqrt{n}}{(\sum_{i=1}^n (Z_i - \overline{Z})^2/(n-1))^{1/2}}, \quad 0 \leq t \leq 1, \tag{2}$$

where the time function $K_n(\cdot)$ is defined as

$$K_n(t) := \sup_{0 \leq m \leq n} \left\{ s_m^2 \leq t s_n^2 \right\}, \quad 0 \leq t \leq 1, \tag{3}$$

with

$$s_0^2 := 0 \quad \text{and} \quad s_m^2 := \sum_{i=1}^m \sigma_i^2, \quad m \geq 1. \tag{4}$$

In (2), we put $\sum_{i=1}^0 Z_i := 0$.

Csörgő, Szyszkowicz and Wang (2003, 2004), among other things, study self-normalized partial sums and a corresponding process in $D[0,1]$ for the above $\{Z_i, i \geq 1\}$, namely

$$V_n(Z_1,\ldots,Z_n) := \frac{\sum_{i=1}^n Z_i}{(\sum_{i=1}^n Z_i^2)^{1/2}} \quad \text{and} \quad V_n^t(Z_1,\ldots,Z_n) := \frac{\sum_{i=1}^{K_n(t)} Z_i}{(\sum_{i=1}^n Z_i^2)^{1/2}}. \tag{5}$$

On assuming that the Lindeberg condition for $Z_i$'s is satisfied, that is

$$\text{for each } \varepsilon > 0, \quad s_n^{-2} \sum_{i=1}^n E Z_i^2 1\!\!1_{\{|Z_i| \geq \varepsilon s_n\}} \to 0, \quad \text{as } n \to \infty, \tag{6}$$

where $1\!\!1_A$ denotes the indicator function of set $A$, as a consequence of a well-known result of Prohorov (1956, Theorem 3.1), they conclude (cf. Proposition 2.2 combined with Remark 2.6 in Csörgő *et al.* (2004)):

$$V_n^t(Z_1,\ldots,Z_n) \xrightarrow{\mathcal{D}} W(t) \quad \text{on} \quad (D[0,1],\rho), \quad n \to \infty, \tag{7}$$

where $\{W(t), 0 \leq t \leq 1\}$ is a standard Wiener process, and $\rho$ stands for the sup-norm metric on $D[0,1]$. The weak convergence in (7) is a weak invariance principle, and it amounts to the following functional central limit theorem (FCLT) (cf., e.g., Sections 3.3 and 3.4 in Csörgő (2002)):

$$h\left(V_n^t(Z_1,\ldots,Z_n)\right) \xrightarrow{\mathcal{D}} h(W(t)), \quad n \to \infty, \tag{8}$$

for all functionals $h : D[0,1] \to \mathbb{R}$ that are $\mathcal{D}$-measurable and $\rho$-continuous, or $\rho$-continuous except at points forming a set of Wiener measure zero on $(D[0,1], \mathcal{D})$, where $\mathcal{D}$ is the sigma-field of subsets of $D[0,1]$ generated by the finite-dimensional subsets of $D[0,1]$.

Csörgő *et al.* (2004) also show (cf. their Proposition 2.3) that one can redefine mean zero $\{Z_i,\, i \geq 1\}$ as in (6) on a richer probability space together with a sequence of independent standard normal r.v.'s $\{Y_i,\, i \geq 1\}$, such that

$$\sup_{0 \leq t \leq 1} \left| V_n^t(Z_1,\ldots,Z_n) - s_n^{-1} \sum_{i=1}^{K_n(t)} \sigma_i Y_i \right| = o_P(1), \quad n \to \infty. \tag{9}$$

Since $s_n^{-1} \sum_{i=1}^{K_n(t)} \sigma_i Y_i \stackrel{\mathcal{D}}{=} W(s_n^{-2} \sum_{i=1}^{K_n(t)} \sigma_i^2)$, $0 \leq t \leq 1$, and, on account of (6), $\sup_{0 \leq t \leq 1} |s_n^{-2} \sum_{i=1}^{K_n(t)} \sigma_i^2 - t| \leq s_n^{-2} \max_{1 \leq i \leq n} \sigma_i^2 \to 0$, $n \to \infty$, then by using the Lévy modulus of continuity of a Wiener process (cf., e.g., Csörgő and Révész (1981)), $\sup_{0 \leq t \leq 1} |W(s_n^{-2} \sum_{i=1}^{K_n(t)} \sigma_i^2) - W(t)| = o_P(1)$, $n \to \infty$. The latter nearness combined with the notion of an FCLT on $(D[0,1],\rho)$ results in

$$s_n^{-1} \sum_{i=1}^{K_n(t)} \sigma_i Y_i \xrightarrow{\mathcal{D}} W(t) \quad \text{on } (D[0,1],\rho), \quad n \to \infty. \tag{10}$$

Consequently, in view of (10), the sup-norm approximation in probability in (9) implies (8), that is the FCLT in (7). Moreover, although the main foci of this paper are FCLT's and their application to constructing asymptotic confidence intervals, upcoming statements of the FCLT's in Lemma 1 and our main Theorem 1 will be accompanied by corresponding sup-norm approximations in probability in (c) parts of Lemma 1 and Theorem 1, as by proving these approximations and using (10), we also establish the FCLT's *à la* (8).

The results in (7) and (9) almost immediately yield Lemma 1 with the corresponding analogues for the Student process $T_n^t(Z_1,\ldots,Z_n)$ of (2).

**Lemma 1.** *Let* $\{Z_i,\, i\geq 1\}$ *be independent mean zero r.v.'s with finite positive variances* $\mathrm{Var} Z_i=\sigma_i^2$, $i\geq 1$. *Assume also the Lindeberg condition as in* (6). *Then, as* $n\to\infty$,

(a) $T_n^{t_0}(Z_1,\ldots,Z_n)\stackrel{\mathcal{D}}{\to} N(0,t_0),\quad t_0\in(0,1]$;

(b) $T_n^t(Z_1,\ldots,Z_n)\stackrel{\mathcal{D}}{\to} W(t)\quad$ *on* $\quad (D[0,1],\rho)$;

(c) *we can redefine* $\{Z_i,\, i\geq 1\}$ *on a richer probability space together with a sequence of independent standard normal r.v.'s* $\{Y_i,\, i\geq 1\}$ *such that*

$$\sup_{0\leq t\leq 1}\left|T_n^t(Z_1,\ldots,Z_n)-s_n^{-1}\sum_{i=1}^{K_n(t)}\sigma_iY_i\right|=o_P(1),$$

*where* $K_n(t)$ *and* $s_n$ *are as in* (3) *and* (4).

*Case of symmetric random variables*

Consider now independent symmetric mean zero r.v.'s $\{Z_i, i\geq 1\}$ that do not necessarily have finite variances. For such r.v.'s, Egorov (1996) proves that, as $n\to\infty$,

$$V_n(Z_1,\ldots,Z_n)\stackrel{\mathcal{D}}{\to} N(0,1)\quad\text{if and only if}\quad \frac{\max_{1\leq i\leq n} Z_i^2}{\sum_{i=1}^n Z_i^2}\stackrel{P}{\to}0. \tag{11}$$

Aiming at a generalization of (11) for an appropriate $D[0,1]$ version of $V_n(Z_1,\ldots, Z_n)$, Csörgő *et al.* (2003) introduce the self-normalized partial sums process

$$\widehat{V}_n^t(Z_1,\ldots,Z_n):=\frac{\sum_{i=1}^{\widehat{K}_n(t)}Z_i}{(\sum_{i=1}^n Z_i^2)^{1/2}},\quad 0\leq t\leq 1, \tag{12}$$

where the time function $\widehat{K}_n(t)$ is a suitable analogue of $K_n(t)$ of (3) for the r.v.'s $Z_i$'s with not necessarily finite variances, namely

$$\widehat{K}_n(t):=\sup_{0\leq m\leq n}\left\{\sum_{i=1}^m Z_i^2\leq t\sum_{i=1}^n Z_i^2\right\},\quad 0\leq t\leq 1,\quad \sum_{i=1}^0 Z_i^2:=0. \tag{13}$$

Csörgő *et al.* (2003, Theorem 2) show that, as $n \to \infty$,

$$\widehat{V}_n^t(Z_1, \ldots, Z_n) \xrightarrow{\mathcal{D}} W(t) \quad \text{on } (D[0,1], \rho) \quad \text{if and only if} \quad \frac{\max_{1 \le i \le n} Z_i^2}{\sum_{i=1}^n Z_i^2} \xrightarrow{P} 0. \tag{14}$$

In view of $\widehat{V}_n^t(Z_1, \ldots, Z_n)$, one can define and study the Student process

$$\widehat{T}_n^t(Z_1, \ldots, Z_n) := \frac{\sum_{i=1}^{\widehat{K}_n(t)} Z_i / \sqrt{n}}{\left(\sum_{i=1}^n (Z_i - \overline{Z})^2/(n-1)\right)^{1/2}}, \quad 0 \le t \le 1, \tag{15}$$

with $\widehat{K}_n(t)$ of (13). It is not hard to see that

$$\widehat{T}_n^t(Z_1, \ldots, Z_n) = \frac{\widehat{V}_n^t(Z_1, \ldots, Z_n)}{\sqrt{\left(n - V_n^2(Z_1, \ldots, Z_n)\right)/(n-1)}}, \quad 0 \le t \le 1. \tag{16}$$

Hence, if $\widehat{T}_n^t(Z_1, \ldots, Z_n)$ or $\widehat{V}_n^t(Z_1, \ldots, Z_n)$ has an asymptotic distribution, then so does the other, and these distributions coincide. Consequently, (14) also holds true for $\widehat{T}_n^t(Z_1, \ldots, Z_n)$.

**Lemma 2.** *Let $\{Z_i, i \ge 1\}$ be independent mean zero symmetric r.v.'s. Then, $\widehat{T}_n^t(Z_1, \ldots, Z_n) \xrightarrow{\mathcal{D}} W(t)$ on $(D[0,1], \rho)$ if and only if $\max_{1 \le i \le n} Z_i^2 / \sum_{i=1}^n Z_i^2 \xrightarrow{P} 0$, as $n \to \infty$.*

### *1.2. Main results: functional asymptotic confidence intervals for a common mean of independent random variables*

*Case of random variables satisfying Lindeberg's condition*

Consider independent r.v.'s $\{Z_i, i \ge 1\}$ with a common mean $\mu$ and finite positive variances that satisfy Lindeberg's condition. As a consequence of the (a) part with $t_0 = 1$ of Lemma 1, $T_n(Z_1 - \mu, \ldots, Z_n - \mu) \xrightarrow{\mathcal{D}} N(0,1)$, $n \to \infty$. Since the Student statistic $T_n(Z_1 - \mu, \ldots, Z_n - \mu)$ does not contain the typically unknown variances $\sigma_i^2$, as compared to the expression $s_n^{-1} \sum_{i=1}^n (Z_i - \mu)$ in the usual statement of the classical Lindeberg-Feller central limit theorem (CLT), the above CLT for $T_n(Z_1 - \mu, \ldots, Z_n - \mu)$ can be used for asymptotic confidence interval (CI) estimation of the mean $\mu$.

The data-based Studentized FCLT in the (b) part of Lemma 1 provides a source of further asymptotic CI's, or what we call functional asymptotic CI's (FACI's), for $\mu$. For example, since the sup-functional $\sup_{0 \le t \le 1} |\cdot|$ on $D[0,1]$ is $\rho$-continuous, from the (b) part of Lemma 1 we conclude

$$\sup_{0 \le t \le 1} \left| T_n^t(Z_1 - \mu, \ldots, Z_n - \mu) \right| \xrightarrow{\mathcal{D}} \sup_{0 \le t \le 1} |W(t)|, \quad n \to \infty, \tag{17}$$

and the latter convergence in distribution yields a $1-\alpha$ size FACI for $\mu$ as follows:

$$\bigcap_{k=1}^{n}\left[\frac{\sum_{i=1}^{k} Z_i - a\sqrt{\frac{n\sum_{i=1}^{n}(Z_i-\overline{Z})^2}{n-1}}}{k}, \frac{\sum_{i=1}^{k} Z_i + a\sqrt{\frac{n\sum_{i=1}^{n}(Z_i-\overline{Z})^2}{n-1}}}{k}\right], \tag{18}$$

where $P\left(\sup_{0\le t\le 1}|W(t)|>a\right)=\alpha$, $0<\alpha<1$. The distribution function of the r.v. $\sup_{0\le t\le 1}|W(t)|$ can be found, for example, in Csörgő and Révész (1981), as well as in Csörgő and Horváth (1984) where it is also tabulated.

In construction of the FACI (18) for $\mu$, due to the very nature of the sup-functional $\sup_{0\le t\le 1}|\cdot|$ on $D[0,1]$, the time function $K_n(t)$ of (3) of $T_n^t(Z_1-\mu,\dots,Z_n-\mu)$ was employed only to the extent of using its values $0,1,\dots,n$. However, when dealing with some other appropriate functionals in regard of constructing FACI's for $\mu$ from the FCLT in (b) of Lemma 1, the jump points $s_k^2/s_n^2$ of the step-function $K_n(t)$ may also enter the picture. These jump points are typically unknown, unless $Z_i$'s have equal variances and $s_k^2/s_n^2=k/n$. To resolve this problem, we replace $K_n(t)$ with its "empirical", data-based version

$$\widetilde{K}_n(t) := \sup_{0\le m\le n}\left\{\sum_{i=1}^{m}(Z_i-\overline{Z})^2 \le t\sum_{i=1}^{n}(Z_i-\overline{Z})^2\right\}, \quad 0\le t\le 1, \tag{19}$$

where $\sum_{i=1}^{0}(Z_i-\overline{Z})^2 := 0$, and establish our main Theorem 1, an analogue of Lemma 1 for the Student process

$$\widetilde{T}_n^t(Z_1,\dots,Z_n) := \frac{\sum_{i=1}^{\widetilde{K}_n(t)} Z_i/\sqrt{n}}{\left(\sum_{i=1}^{n}(Z_i-\overline{Z})^2/(n-1)\right)^{1/2}}, \quad 0\le t\le 1. \tag{20}$$

Without loss of generality, Theorem 1 is stated under the assumption that $\mu=0$.

**Theorem 1.** *Let $\{Z_i,\ i\ge 1\}$ be independent mean zero r.v.'s with finite positive variances $\mathrm{Var}Z_i=\sigma_i^2$, $i\ge 1$. Assume also the Lindeberg condition as in (6). Then, as $n\to\infty$,*

(a) $\widetilde{T}_n^{t_0}(Z_1,\dots,Z_n) \xrightarrow{\mathcal{D}} N(0,t_0), \quad t_0\in(0,1];$
(b) $\widetilde{T}_n^{t}(Z_1,\dots,Z_n) \xrightarrow{\mathcal{D}} W(t)$ *on* $(D[0,1],\rho);$
(c) *we can redefine $\{Z_i,\ i\ge 1\}$ on a richer probability space together with a sequence of independent standard normal r.v.'s $\{Y_i,\ i\ge 1\}$ such that*

$$\sup_{0\le t\le 1}\left|\widetilde{T}_n^t(Z_1,\dots,Z_n) - s_n^{-1}\sum_{i=1}^{K_n(t)}\sigma_i Y_i\right| = o_P(1),$$

*where $K_n(t)$ and $s_n$ are as in (3) and (4).*

To illustrate when construction of FACI's for a not necessarily zero $\mu$ call for the FCLT as in the (b) part of Theorem 1, we consider convergence in distribution of two special functionals of $\widetilde{T}_n^t(Z_1-\mu,\dots,Z_n-\mu)$ in Examples 1 and 2.

**Example 1.** For a fixed $t_0 \in (0,1]$, we consider $\widetilde{T}_n^{t_0}(Z_1-\mu,\ldots,Z_n-\mu)$, one of the simplest $\rho$-continuous functionals of $\widetilde{T}_n^{t}(Z_1-\mu,\ldots,Z_n-\mu)$. As a consequence of the FCLT in (b) of Theorem 1, or directly by (a) of Theorem 1, we obtain the following $1-\alpha$ size FACI for $\mu$:

$$\left[\frac{\sum_{i=1}^{\widetilde{K}_n(t_0)} Z_i - z_{\alpha/2}\sqrt{t_0}\sqrt{\frac{n\sum_{i=1}^n (Z_i-\overline{Z})^2}{n-1}}}{\widetilde{K}_n(t_0)}, \frac{\sum_{i=1}^{\widetilde{K}_n(t_0)} Z_i + z_{\alpha/2}\sqrt{t_0}\sqrt{\frac{n\sum_{i=1}^n (Z_i-\overline{Z})^2}{n-1}}}{\widetilde{K}_n(t_0)}\right], \tag{21}$$

where $z_{\alpha/2}$ is the $100(1-\alpha/2)^{\text{th}}$ percentile of the standard normal distribution. The FACI in (21) is completely data-based, as $\widetilde{K}_n(t_0)$ is computable. Indeed, if $t_0 = 1$, then $\widetilde{K}_n(t_0) = n$, while for $t_0 \in (0,1)$ and a given sample $Z_1,\ldots,Z_n$, we can find $k_0$, $0 \le k_0 \le n-1$, such that $\sum_{i=1}^{k_0}(Z_i-\overline{Z})^2/\sum_{i=1}^n(Z_i-\overline{Z})^2 \le t_0 < \sum_{i=1}^{k_0+1}(Z_i-\overline{Z})^2/\sum_{i=1}^n(Z_i-\overline{Z})^2$, and consequently, $\widetilde{K}_n(t_0) = k_0$. We also note that (21) is well-defined, as by (42) below, $\max_{1\le i\le n}(Z_i-\overline{Z})^2/\sum_{i=1}^n(Z_i-\overline{Z})^2 \stackrel{P}{\to} 0$ and hence, $(Z_1-\overline{Z})^2/\sum_{i=1}^n(Z_i-\overline{Z})^2 < t_0$ (or $\widetilde{K}_n(t_0) \ne 0$) with probability approaching one, as $n\to\infty$.

**Example 2.** The integral functional $\int_0^1 \cdot\, dt$ on $D[0,1]$ is $\rho$-continuous, as for any $f(t)$ and $g(t)$ in $D[0,1]$, $|\int_0^1 f(t)dt - \int_0^1 g(t)dt| \le \sup_{0\le t\le 1}|f(t)-g(t)|$. In view of this and the FCLT in (b) part of Theorem 1, as $n\to\infty$,

$$\int_0^1 \widetilde{T}_n^t(Z_1-\mu,\ldots,Z_n-\mu)dt \stackrel{\mathcal{D}}{\to} \int_0^1 W(t)dt \stackrel{\mathcal{D}}{=} N(0,1/3). \tag{22}$$

By noting that

$$\int_0^1 \widetilde{T}_n^t(Z_1-\mu,\ldots,Z_n-\mu)dt = \sum_{k=1}^{n-1}\nu_{k+1}\frac{\sum_{i=1}^k (Z_i-\mu)/\sqrt{n}}{(\sum_{i=1}^n (Z_i-\overline{Z})^2/(n-1))^{1/2}},$$

with

$$\nu_k := \frac{(Z_k-\overline{Z})^2}{\sum_{i=1}^n (Z_i-\overline{Z})^2}, \quad 1\le k\le n, \tag{23}$$

we obtain a $1-\alpha$ size FACI for $\mu$ with the lower and upper bounds given by

$$\frac{\sum_{k=1}^{n-1}\nu_{k+1}\sum_{i=1}^k Z_i \mp \frac{z_{\alpha/2}}{\sqrt{3}}\sqrt{\frac{n\sum_{i=1}^n (Z_i-\overline{Z})^2}{n-1}}}{\sum_{k=1}^{n-1}\nu_{k+1}k}, \tag{24}$$

where $P\left(|N(0,1/3)| > z_{\alpha/2}/\sqrt{3}\right) = \alpha$.

It would naturally be desirable to investigate individual and comparative performances, such as the expected lengths for example, of the obtained FACI's for $\mu$ in (18), (21) and (24).

*Case of symmetric random variables*

Consider independent symmetric r.v.'s $Z_i$ with a common, not necessarily zero, mean $\mu$. Lemma 2 remains true for such $Z_i$'s if $Z_1, \ldots, Z_n$ are replaced with $Z_1 - \mu, \ldots, Z_n - \mu$ in its statement. However, in the thus stated Lemma 2, the time function of $\widehat{T}_n^t(Z_1 - \mu, \ldots, Z_n - \mu)$ becomes $\sup_{0 \le m \le n}\{\sum_{i=1}^m (Z_i - \mu)^2 \le t \sum_{i=1}^n (Z_i - \mu)^2\}$, a function of an unknown $\mu$. Hence, such an FCLT is not necessarily of immediate use for construction of various FACI's for $\mu$, just like the FCLT of (b) of Lemma 1 when it is applied to independent mean $\mu$ r.v.'s $Z_i$ satisfying Lindeberg's condition (cf. the lines preceding (19)). To extend the applicability of the FCLT of Lemma 2 in this regard, we first establish our main Theorem 2 for the Student process $\widetilde{T}_n^t(Z_1 - \mu, \ldots, Z_n - \mu)$ as in (20), a data-based version of Lemma 2 that uses $\widetilde{K}_n(t)$ of (19) instead of the above noted time function of $\widehat{T}_n^t(Z_1 - \mu, \ldots, Z_n - \mu)$.

**Theorem 2.** *Let $\{Z_i, i \ge 1\}$ be independent symmetric r.v.'s with a common mean $\mu$. Then, for $\widetilde{T}_n^t(Z_1, \ldots, Z_n)$ as in (20), $\widetilde{T}_n^t(Z_1 - \mu, \ldots, Z_n - \mu) \xrightarrow{\mathcal{D}} W(t)$ on $(D[0,1], \rho)$ if and only if $\max_{1 \le i \le n}(Z_i - \mu)^2 / \sum_{i=1}^n (Z_i - \mu)^2 \xrightarrow{P} 0$, as $n \to \infty$.*

It is appealing to replace the convergence $\max_{1 \le i \le n}(Z_i - \mu)^2 / \sum_{i=1}^n (Z_i - \mu)^2 \xrightarrow{P} 0$ in Theorem 2 with the data-based one of $\max_{1 \le i \le n}(Z_i - \overline{Z})^2 / \sum_{i=1}^n (Z_i - \overline{Z})^2 \xrightarrow{P} 0$, as $n \to \infty$, especially when concluding FACI's for $\mu$ via appropriate functionals of $\widetilde{T}_n^t(Z_1 - \mu, \ldots, Z_n - \mu)$ of Theorem 2. Hence we present Corollary 1 that amounts to a completely data-based version of Theorem 2.

**Corollary 1.** *Let $\{Z_i, i \ge 1\}$ be independent symmetric r.v.'s with a common mean $\mu$. Then $\widetilde{T}_n^t(Z_1 - \mu, \ldots, Z_n - \mu) \xrightarrow{\mathcal{D}} W(t)$ on $(D[0,1], \rho)$ if and only if $\max_{1 \le i \le n}(Z_i - \overline{Z})^2 / \sum_{i=1}^n (Z_i - \overline{Z})^2 \xrightarrow{P} 0$, as $n \to \infty$.*

We note that, in view of Corollary 1, the FACI's for $\mu$ in (21) and (24) also hold true for independent symmetric r.v.'s $Z_i$ with a common mean $\mu$ and not necessarily finite variances, provided that $\max_{1 \le i \le n}(Z_i - \overline{Z})^2 / \sum_{i=1}^n (Z_i - \overline{Z})^2 \xrightarrow{P} 0$, as $n \to \infty$.

## 2. Proofs

Hereafter, notations $o_P(1)$ and $O_P(1)$ stand for sequences of r.v.'s that, respectively, converge to zero and are bounded in probability, as $n \to \infty$.

**Proof of Lemma 1.** In view of (10), the proof reduces to establishing the (c) part of Lemma 1. On account of (7), $\sup_{0 \le t \le 1} |V_n^t(Z_1, \ldots, Z_n)| = O_P(1)$ and $V_n(Z_1, \ldots, Z_n) = O_P(1)$. Moreover, we also have (9) with $\{Z_i, i \ge 1\}$ and $\{Y_i, i \ge 1\}$ defined on the same probability space. Combining all this with a representation for $T_n^t(Z_1, \ldots, Z_n)$ *à la* (16), as $n \to \infty$, we arrive at

$$\sup_{0 \le t \le 1} \left| T_n^t(Z_1, \ldots, Z_n) - \frac{\sum_{i=1}^{K_n(t)} \sigma_i Y_i}{s_n} \right| \le \sup_{0 \le t \le 1} \left| V_n^t(Z_1, \ldots, Z_n) - \frac{\sum_{i=1}^{K_n(t)} \sigma_i Y_i}{s_n} \right|$$

$$+ \sup_{0\le t\le 1} \left|V_n^t(Z_1,\ldots,Z_n)\right| \left| \frac{1}{\sqrt{(n - V_n^2(Z_1,\ldots,Z_n))/(n-1)}} - 1 \right|$$
$$= o_P(1) + O_P(1)o_P(1) = o_P(1). \tag{25}$$

□

Next, we spell out a special case of Raikov's theorem (cf. Theorem 4 on p.143 in Gnedenko and Kolmogorov (1954)), and then, establish auxiliary Lemma 4 that is required for the proof of Theorem 1.

**Lemma 3.** *Let $\{Z_i,\, i \ge 1\}$ be independent mean zero r.v.'s with finite positive variances* $\mathrm{Var} Z_i = \sigma_i^2$, $i \ge 1$. *Suppose also that the Lindeberg condition in (6) is satisfied. Then,*

$$s_n^{-2} \sum_{i=1}^{n} Z_i^2 \xrightarrow{P} 1, \quad n \to \infty, \tag{26}$$

*with $s_n^2$ of (4).*

**Lemma 4.** *For $\{Z_i,\, i \ge 1\}$ as in Lemma 3,*

$$\sup_{0\le t\le 1} \left| V_n^t(Z_1,\ldots,Z_n) - \widetilde{V}_n^t(Z_1,\ldots,Z_n) \right| = o_P(1), \quad n \to \infty, \tag{27}$$

*where the self-normalized partial sums processes $V_n^t(Z_1,\ldots,Z_n)$ and $\widetilde{V}_n^t(Z_1,\ldots,Z_n)$ are defined as follows: for $0 \le t \le 1$,*

$$V_n^t(Z_1,\ldots,Z_n) = \frac{\sum_{i=1}^{K_n(t)} Z_i}{(\sum_{i=1}^n Z_i^2)^{1/2}} \quad \text{and} \quad \widetilde{V}_n^t(Z_1,\ldots,Z_n) = \frac{\sum_{i=1}^{\widetilde{K}_n(t)} Z_i}{(\sum_{i=1}^n Z_i^2)^{1/2}}, \tag{28}$$

*with the time functions $K_n(t)$ and $\widetilde{K}_n(t)$ as in (3) and (19).*

**Proof of Lemma 4.** The scheme of this proof is motivated by the lines of the proof of Theorem 2 in Račkauskas and Suquet (2001).

Let $U_n^t$ and $\widetilde{U}_n^t$ be the respective $C[0,1]$ "Donskerized" versions of the $D[0,1]$ processes $V_n^t(Z_1,\cdots,Z_n)$ and $\widetilde{V}_n^t(Z_1,\cdots,Z_n)$. Namely, $U_n^t$ and $\widetilde{U}_n^t$, as continuous functions of $t$, are linear respectively on the intervals $[s_k^2/s_n^2, s_{k+1}^2/s_n^2]$ and $[\sum_{i=1}^k (Z_i-\overline{Z})^2/\sum_{i=1}^n (Z_i-\overline{Z})^2, \sum_{i=1}^{k+1}(Z_i-\overline{Z})^2/\sum_{i=1}^n (Z_i-\overline{Z})^2]$ for each $k = 0,1,\cdots,n-1$, and both take values $\sum_{i=1}^k Z_i/(\sum_{i=1}^n Z_i^2)^{1/2}$ respectively at $s_k^2/s_n^2$ and $\sum_{i=1}^k (Z_i-\overline{Z})^2/\sum_{i=1}^n (Z_i-\overline{Z})^2$, $k = 0,1,\cdots,n$, where $s_k^2$ is defined in (4), and $\sum_{i=1}^0 (Z_i-\overline{Z})^2 := 0$ and $\sum_{i=1}^0 Z_i := 0$. Clearly,

$$\begin{aligned} \sup_{0\le t\le 1} \left|V_n^t(Z_1,\ldots,Z_n) - U_n^t\right| &\le \frac{\max_{1\le i\le n}|Z_i|}{(\sum_{i=1}^n Z_i^2)^{1/2}}, \\ \sup_{0\le t\le 1} \left|\widetilde{V}_n^t(Z_1,\ldots,Z_n) - \widetilde{U}_n^t\right| &\le \frac{\max_{1\le i\le n}|Z_i|}{(\sum_{i=1}^n Z_i^2)^{1/2}}. \end{aligned} \tag{29}$$

From Lindeberg's condition in (6), for any $\varepsilon > 0$, as $n \to \infty$,

$$P\left(\max_{1\le i\le n} |Z_i| \ge \varepsilon s_n\right) \le \sum_{i=1}^{n} P(|Z_i| \ge \varepsilon s_n) \le (\varepsilon s_n)^{-2}\sum_{i=1}^{n} E\left(Z_i^2 \mathbb{1}_{\{|Z_i|\ge \varepsilon s_n\}}\right) \to 0. \tag{30}$$

On combining (26) and (30),

$$\frac{\max_{1\le i\le n}|Z_i|}{(\sum_{i=1}^n Z_i^2)^{1/2}} = o_P(1), \quad n\to\infty. \tag{31}$$

In view of (29) and (31), in order to show (27), it suffices to prove that

$$\sup_{0\le t\le 1}\left|U_n^t - \widetilde{U}_n^t\right| = o_P(1), \quad n\to\infty. \tag{32}$$

Let $\theta_n(t)$ be the random element of $C[0,1]$ that is linear in $t$ on the intervals $[s_k^2/s_n^2, s_{k+1}^2/s_n^2]$ for each $k = 0, 1, \ldots, n-1$, with $\theta_n(s_k^2/s_n^2) = \sum_{i=1}^k (Z_i - \overline{Z})^2 / \sum_{i=1}^n (Z_i - \overline{Z})^2$, $k = 0, 1, \ldots, n$. Since $\sup_{0\le t\le 1}|U_n^t - \widetilde{U}_n^t| = \sup_{0\le t\le 1}|U_n^{\theta_n(t)} - \widetilde{U}_n^{\theta_n(t)}|$ and $\widetilde{U}_n^{\theta_n(t)} = U_n^t$, $0\le t\le 1$, then (32) reads as

$$\sup_{0\le t\le 1}\left|U_n^{\theta_n(t)} - U_n^t\right| = o_P(1), \quad n\to\infty. \tag{33}$$

For $f(t)\in C[0,1]$, let $\omega(f(t);\delta) := \sup_{|t_1-t_2|\le\delta}|f(t_1)-f(t_2)|$ be the modulus of continuity of $f(t)$. For any $\lambda > 0$ and $0 < \delta \le 1$, we have

$$\begin{aligned}
&P\left(\sup_{0\le t\le 1}|U_n^{\theta_n(t)} - U_n^t| \ge \lambda\right) \\
&\le P\left(\sup_{0\le t\le 1,\, |t-\theta_n(t)|\le \sup_{0\le t\le 1}|t-\theta_n(t)|}|U_n^{\theta_n(t)} - U_n^t| \ge \lambda\right) \\
&\le P\left(\omega(U_n^t;\delta)\ge\lambda\right) + P\left(\sup_{0\le t\le 1}|t-\theta_n(t)| > \delta\right).
\end{aligned} \tag{34}$$

By Theorem 3.1 in Prohorov (1956) and (26), $U_n^t \xrightarrow{\mathcal{D}} W(t)$ on $(C[0,1],\rho)$, $n\to\infty$, and therefore, for the continuous functional $\omega(\cdot;\delta)$ on $C[0,1]$ and any $\lambda>0$,

$$P\left(\omega(U_n^t;\delta)\ge\lambda\right) \to P\left(\omega(W(t);\delta)\ge\lambda\right), \quad n\to\infty. \tag{35}$$

In view of, for example, the Lévy modulus of continuity of a Wiener process (cf., e.g., Csörgő and Révész (1981)), for any $\varepsilon>0$ there is $\delta\in(0,1]$ such that

$$P\left(\omega(W(t);\delta)\ge\lambda\right) < \varepsilon. \tag{36}$$

Taking into account (34)–(36), to complete the proof of Lemma 4, we only need to verify that

$$\sup_{0\le t\le 1}|t-\theta_n(t)| = o_P(1), \quad n\to\infty. \tag{37}$$

To prove (37), we first write

$$\sup_{0\le t\le 1}|t-\theta_n(t)| \le \sup_{0\le t\le 1}\left|t-\frac{\sum_{i=1}^{K_n(t)}\sigma_i^2}{s_n^2}\right| + \sup_{0\le t\le 1}\left|\frac{\sum_{i=1}^{K_n(t)}\sigma_i^2}{s_n^2}-\frac{\sum_{i=1}^{K_n(t)}Z_i^2}{s_n^2}\right|$$
$$+ \sup_{0\le t\le 1}\left|\frac{\sum_{i=1}^{K_n(t)}Z_i^2}{s_n^2}-\frac{\sum_{i=1}^{K_n(t)}(Z_i-\overline{Z})^2}{s_n^2}\right|$$
$$+ \sup_{0\le t\le 1}\left|\frac{\sum_{i=1}^{K_n(t)}(Z_i-\overline{Z})^2}{s_n^2}-\frac{\sum_{i=1}^{K_n(t)}(Z_i-\overline{Z})^2}{\sum_{i=1}^{n}(Z_i-\overline{Z})^2}\right|$$
$$+ \sup_{0\le t\le 1}\left|\frac{\sum_{i=1}^{K_n(t)}(Z_i-\overline{Z})^2}{\sum_{i=1}^{n}(Z_i-\overline{Z})^2}-\theta_n(t)\right| =: A_1+A_2+A_3+A_4+A_5. \qquad (38)$$

Concerning $A_1$, we have

$$A_1 = \sup_{0\le t\le 1}\left|t-\frac{\sum_{i=1}^{K_n(t)}\sigma_i^2}{s_n^2}\right| \le \frac{\max_{1\le i\le n}\sigma_i^2}{s_n^2}\to 0, \quad n\to\infty, \qquad (39)$$

on account of Feller's condition for $Z_i$'s that is implied by Lindeberg's one in (6). From (7) and (26), $|\sum_{i=1}^{n}Z_i|/s_n = O_P(1)$ and $\sup_{0\le t\le 1}|\sum_{i=1}^{K_n(t)}Z_i|/s_n = O_P(1)$. Using these facts together with (26) and (31), we have that for $A_3$, $A_4$ and $A_5$ in (38), as $n\to\infty$,

$$A_3 = \sup_{0\le t\le 1}\left|\frac{\sum_{i=1}^{K_n(t)}Z_i^2}{s_n^2}-\frac{\sum_{i=1}^{K_n(t)}(Z_i-\overline{Z})^2}{s_n^2}\right| \le \frac{2|\sum_{i=1}^{n}Z_i|}{ns_n}\sup_{0\le t\le 1}\frac{|\sum_{i=1}^{K_n(t)}Z_i|}{s_n}$$
$$+\frac{1}{n}\left(\frac{\sum_{i=1}^{n}Z_i}{s_n}\right)^2 = \frac{O_P(1)}{n}O_P(1)+\frac{O_P(1)}{n} = o_P(1), \qquad (40)$$

$$A_4 = \sup_{0\le t\le 1}\left|\frac{\sum_{i=1}^{K_n(t)}(Z_i-\overline{Z})^2}{s_n^2}-\frac{\sum_{i=1}^{K_n(t)}(Z_i-\overline{Z})^2}{\sum_{i=1}^{n}(Z_i-\overline{Z})^2}\right| \le \left|\frac{\sum_{i=1}^{n}(Z_i-\overline{Z})^2}{s_n^2}-1\right|$$
$$\le \left|\frac{\sum_{i=1}^{n}Z_i^2}{s_n^2}-1\right|+\frac{1}{n}\left(\frac{\sum_{i=1}^{n}Z_i}{s_n}\right)^2 = o_P(1)+\frac{O_P(1)}{n} = o_P(1), \qquad (41)$$

$$A_5 = \sup_{0\le t\le 1}\left|\frac{\sum_{i=1}^{K_n(t)}(Z_i-\overline{Z})^2}{\sum_{i=1}^{n}(Z_i-\overline{Z})^2}-\theta_n(t)\right| \le \frac{\max_{1\le i\le n}(Z_i-\overline{Z})^2}{\sum_{i=1}^{n}(Z_i-\overline{Z})^2}$$
$$\le \frac{\sum_{i=1}^{n}Z_i^2}{s_n^2}\,\frac{s_n^2}{\sum_{i=1}^{n}(Z_i-\overline{Z})^2}\,\frac{4\max_{1\le i\le n}Z_i^2}{\sum_{i=1}^{n}Z_i^2}$$
$$= (1+o_P(1))(1+o_P(1))o_P(1) = o_P(1). \qquad (42)$$

Now, in view of (38)–(42), it remains to show that

$$A_2 = \sup_{0\le t\le 1}\left|\frac{\sum_{i=1}^{K_n(t)}\sigma_i^2}{s_n^2}-\frac{\sum_{i=1}^{K_n(t)}Z_i^2}{s_n^2}\right| = o_P(1), \quad n\to\infty. \qquad (43)$$

Let

$$\widetilde{Z}_i^2 = s_n^{-2} Z_i^2 \mathbb{1}_{\{Z_i^2 \le s_n^2\}}, \tag{44}$$

then by Lindeberg's condition in (6) and the inequality

$$\begin{aligned} &\sup_{0\le t\le 1}\left|\frac{\sum_{i=1}^{K_n(t)}\sigma_i^2}{s_n^2} - \frac{\sum_{i=1}^{K_n(t)} Z_i^2}{s_n^2}\right| \\ &\le \sup_{0\le t\le 1}\left|\sum_{i=1}^{K_n(t)}(\widetilde{Z}_i^2 - E\,\widetilde{Z}_i^2)\right| + s_n^{-2}\sum_{i=1}^{n} Z_i^2 \mathbb{1}_{\{Z_i^2 > s_n^2\}} + s_n^{-2}\sum_{i=1}^{n} E\Big(Z_i^2 \mathbb{1}_{\{Z_i^2 > s_n^2\}}\Big), \end{aligned} \tag{45}$$

the proof of (43) narrows down to establishing

$$\sup_{0\le t\le 1}\left|\sum_{i=1}^{K_n(t)}(\widetilde{Z}_i^2 - E\,\widetilde{Z}_i^2)\right| = o_P(1), \quad n\to\infty. \tag{46}$$

By the Ottavani-Skorohod inequality (cf., e.g., Shorack (2000)), for any $a > 0$,

$$\begin{aligned} &P\left(\sup_{0\le t\le 1}\left|\sum_{i=1}^{K_n(t)}(\widetilde{Z}_i^2 - E\,\widetilde{Z}_i^2)\right| \ge 4a\right) = P\left(\max_{1\le k\le n}\left|\sum_{i=1}^{k}(\widetilde{Z}_i^2 - E\,\widetilde{Z}_i^2)\right| \ge 4a\right) \\ &\le \frac{P\Big(\Big|\sum_{i=1}^{n}(\widetilde{Z}_i^2 - E\,\widetilde{Z}_i^2)\Big| \ge 2a\Big)}{1 - \max_{1\le k\le n} P\Big(\Big|\sum_{i=1}^{n}(\widetilde{Z}_i^2 - E\,\widetilde{Z}_i^2) - \sum_{i=1}^{k}(\widetilde{Z}_i^2 - E\,\widetilde{Z}_i^2)\Big| > 2a\Big)}. \end{aligned} \tag{47}$$

On account of the Linderberg condition in (6) and (26),

$$\begin{aligned} \left|\sum_{i=1}^{n}(\widetilde{Z}_i^2 - E\,\widetilde{Z}_i^2)\right| &\le \left|1 - \frac{\sum_{i=1}^{n} Z_i^2}{s_n^2}\right| + s_n^{-2}\sum_{i=1}^{n} Z_i^2 \mathbb{1}_{\{Z_i^2 > s_n^2\}} \\ &+ s_n^{-2}\sum_{i=1}^{n} E\Big(Z_i^2 \mathbb{1}_{\{Z_i^2 > s_n^2\}}\Big) = o_P(1). \end{aligned} \tag{48}$$

If, additionally, for any $a > 0$,

$$\max_{1\le k\le n} P\left(\left|\sum_{i=1}^{k}(\widetilde{Z}_i^2 - E\,\widetilde{Z}_i^2)\right| > a\right) \le \frac{1}{2} \quad \text{for } n \text{ large enough}, \tag{49}$$

then (46) follows from (47)–(49). For $\tau \in (0,1)$ and any $a > 0$,

$$\begin{aligned} &\max_{1\le k\le n} P\left(\left|\sum_{i=1}^{k}(\widetilde{Z}_i^2 - E\,\widetilde{Z}_i^2)\right| > a\right) \le a^{-2}\max_{1\le k\le n} E\left|\sum_{i=1}^{k}(\widetilde{Z}_i^2 - E\,\widetilde{Z}_i^2)\right|^2 \\ &= a^{-2}\sum_{i=1}^{n} E(\widetilde{Z}_i^2 - E\widetilde{Z}_i^2)^2 \le a^{-2}\sum_{i=1}^{n} E\widetilde{Z}_i^4 \end{aligned}$$

$$= a^{-2}s_n^{-4}\sum_{i=1}^{n} E\left(Z_i^4 1\!\!1_{\{Z_i^2<\tau^2 s_n^2\}\cup\{\tau^2 s_n^2\le Z_i^2\le s_n^2\}}\right)$$
$$\le a^{-2}\tau^2 s_n^{-2}\sum_{i=1}^{n} E\left(Z_i^2 1\!\!1_{\{Z_i^2<\tau^2 s_n^2\}}\right) + a^{-2}s_n^{-2}\sum_{i=1}^{n} E\left(Z_i^2 1\!\!1_{\{Z_i^2\ge\tau^2 s_n^2\}}\right)$$
$$\le 2a^{-2}\tau^2 + \frac{1}{4}, \tag{50}$$

where the latter inequality holds for $n$ large enough (depending on $a$ and $\tau$) on account of Lindeberg's condition in (6). In (50), if $a \ge 2\sqrt{2}$, then $2a^{-2}\tau^2+1/4 \le \tau^2/4 + 1/4 < 1/2$ for any $\tau \in (0,1)$, while for $0 < a < 2\sqrt{2}$, by choosing $\tau = a/2\sqrt{2}$, $2a^{-2}\tau^2 + 1/4 = 1/2$. Thus, for any $a > 0$, we achieve (49). This also completes the proof of Lemma 4. □

**Proof of Theorem 1.** Having (10), we are only concerned with the proof of the (c) part of Theorem 1. Representations for $T_n^t(Z_1,\ldots,Z_n)$ and $\widetilde{T}_n^t(Z_1,\ldots,Z_n)$ that rhyme with (16), the (c) part of Lemma 1, (27) of Lemma 4 and (7) result in, as $n \to \infty$,

$$\sup_{0\le t\le 1}\left|\widetilde{T}_n^t(Z_1,\ldots,Z_n) - s_n^{-1}\sum_{i=1}^{K_n(t)}\sigma_i Y_i\right| \le \sup_{0\le t\le 1}\left|T_n^t(Z_1,\ldots,Z_n) - s_n^{-1}\sum_{i=1}^{K_n(t)}\sigma_i Y_i\right|$$
$$+\frac{\sup_{0\le t\le 1}\left|V_n^t(Z_1,\ldots,Z_n) - \widetilde{V}_n^t(Z_1,\ldots,Z_n)\right|}{\sqrt{\left(n - V_n^2(Z_1,\ldots,Z_n)\right)/(n-1)}} = o_P(1) + O_P(1)o_P(1) = o_P(1).$$

□

**Proof of Theorem 2.** Without loss of generality, we assume that $\mu = 0$.

In view of (11) and a version of (16) for $\widetilde{T}_n^t(Z_1,\ldots,Z_n)$, it suffices to prove that $\widetilde{V}_n^t(Z_1,\ldots,Z_n) \xrightarrow{\mathcal{D}} W(t)$ on $(D[0,1],\rho)$ for $\widetilde{V}_n^t(Z_1,\ldots,Z_n)$ of (28), assuming that $\max_{1\le i\le n} Z_i^2/\sum_{i=1}^n Z_i^2 \xrightarrow{P} 0$, as $n \to \infty$. The latter convergence in distribution reduces to establishing that

$$\sup_{0\le t\le 1}\left|\widetilde{V}_n^t(Z_1,\ldots,Z_n) - \widehat{V}_n^t(Z_1,\ldots,Z_n)\right| = o_P(1), \quad n\to\infty, \tag{51}$$

where $\widehat{V}_n^t(Z_1,\ldots,Z_n)$ of (12) satisfies (14). By arguments similar to those used for proving (27) of Lemma 4, (51) follows from

$$\sup_{0\le t\le 1}\left|\widetilde{U}_n^t - \widehat{U}_n^t\right| = o_P(1), \quad n\to\infty, \tag{52}$$

where $\widetilde{U}_n^t$ and $\widehat{U}_n^t$ are the respective $C[0,1]$ ("Donskerized") versions of the $D[0,1]$ processes $\widetilde{V}_n^t(Z_1,\ldots,Z_n)$ and $\widehat{V}_n^t(Z_1,\ldots,Z_n)$.

Consider the random broken line $\eta_n(t) \in C[0,1]$ that is linear in $t$ on the intervals $[\sum_{i=1}^k Z_i^2/\sum_{i=1}^n Z_i^2, \sum_{i=1}^{k+1} Z_i^2/\sum_{i=1}^n Z_i^2]$ for each $k = 0,1,\ldots,n-1$, with $\eta_n(\sum_{i=1}^k Z_i^2/\sum_{i=1}^n Z_i^2)=\sum_{i=1}^k (Z_i-\overline{Z})^2/\sum_{i=1}^n (Z_i-\overline{Z})^2$, $k=0,1,\ldots,n$. Since $\widetilde{U}_n^{\eta_n(t)}=\widehat{U}_n^t$, $0\le t\le 1$, then (52) is equivalent to

$$\sup_{0\le t\le 1}\left|\widehat{U}_n^t-\widehat{U}_n^{\eta_n(t)}\right|=o_P(1),\quad n\to\infty. \tag{53}$$

Via lines similar to those in (34)–(36) and the weak convergence $\widehat{U}_n^t \stackrel{\mathcal{D}}{\to} W(t)$ on $(C[0,1],\rho)$ that follows from the slightly modified proof of the FCLT in (14) for $\widehat{V}_n^t(Z_1,\ldots,Z_n)$ in Csörgő *et al.* (2003), the nearness in (53) is hinged on showing that

$$\sup_{0\le t\le 1}|t-\eta_n(t)|=o_P(1),\quad n\to\infty. \tag{54}$$

We have

$$\begin{aligned}\sup_{0\le t\le 1}|t-\eta_n(t)| &\le \sup_{0\le t\le 1}\left|t-\frac{\sum_{i=1}^{\widehat{K}_n(t)}(Z_i-\overline{Z})^2}{\sum_{i=1}^n (Z_i-\overline{Z})^2}\right|+\frac{\max_{1\le i\le n}(Z_i-\overline{Z})^2}{\sum_{i=1}^n (Z_i-\overline{Z})^2}\\ &=: B_1+B_2, \end{aligned}\tag{55}$$

where, due to (11), (31) and (14),

$$\begin{aligned}B_2 &\le \frac{4\max_{1\le i\le n}Z_i^2}{\sum_{i=1}^n Z_i^2}\frac{\sum_{i=1}^n Z_i^2}{\sum_{i=1}^n (Z_i-\overline{Z})^2}=\frac{o_P(1)}{1-V_n^2(Z_1,\ldots,Z_n)/n}\\ &= \frac{o_P(1)}{1+O_P(1)/n}=o_P(1),\end{aligned}\tag{56}$$

and

$$\begin{aligned}B_1 &\le \sup_{0\le t\le 1}\left|t-\frac{\sum_{i=1}^{\widehat{K}_n(t)}Z_i^2}{\sum_{i=1}^n Z_i^2}\right|+\sup_{0\le t\le 1}\left|\frac{\sum_{i=1}^{\widehat{K}_n(t)}Z_i^2}{\sum_{i=1}^n Z_i^2}-\frac{\sum_{i=1}^{\widehat{K}_n(t)}(Z_i-\overline{Z})^2}{\sum_{i=1}^n Z_i^2}\right|\\ &\quad+\sup_{0\le t\le 1}\left|\frac{\sum_{i=1}^{\widehat{K}_n(t)}(Z_i-\overline{Z})^2}{\sum_{i=1}^n Z_i^2}-\frac{\sum_{i=1}^{\widehat{K}_n(t)}(Z_i-\overline{Z})^2}{\sum_{i=1}^n (Z_i-\overline{Z})^2}\right|\le\frac{\max_{1\le i\le n}Z_i^2}{\sum_{i=1}^n Z_i^2}\\ &\quad+\left(\frac{2|V_n(Z_1,\ldots,Z_n)|}{n}\sup_{0\le t\le 1}\left|\widehat{V}_n^t(Z_1,\ldots,Z_n)\right|+\frac{V_n^2(Z_1,\ldots,Z_n)}{n}\right)\\ &\quad+\left|\frac{\sum_{i=1}^n (Z_i-\overline{Z})^2}{\sum_{i=1}^n Z_i^2}-1\right|\\ &= o_P(1)+\left(\frac{O_P(1)}{n}O_P(1)+\frac{O_P(1)}{n}\right)+\frac{O_P(1)}{n}=o_P(1).\end{aligned}\tag{57}$$

□

**Proof of Corollary 1.** Let $\widetilde{T}_n^t(Z_1-\mu,\ldots,Z_n-\mu)\stackrel{\mathcal{D}}{\to}W(t)$ on $(D[0,1],\rho)$, $n\to\infty$. In view of Theorem 2 and a representation for $\widetilde{T}_n^t(Z_1-\mu,\ldots,Z_n-\mu)$

*à la* (16), this FCLT yields $\max_{1\le i\le n}(Z_i-\mu)^2/\sum_{i=1}^n(Z_i-\mu)^2 = o_P(1)$ and $\sum_{i=1}^n(Z_i-\mu)/(\sum_{i=1}^n(Z_i-\mu)^2)^{1/2} = O_P(1)$, $n\to\infty$. Consequently,

$$\frac{\max_{1\le i\le n}(Z_i-\overline{Z})^2}{\sum_{i=1}^n(Z_i-\overline{Z})^2} \le \frac{2\max_{1\le i\le n}(Z_i-\mu)^2+2(\overline{Z}-\mu)^2}{\sum_{i=1}^n(Z_i-\mu)^2-n(\overline{Z}-\mu)^2}$$

$$= \frac{2\max\limits_{1\le i\le n}(Z_i-\mu)^2/\sum_{i=1}^n(Z_i-\mu)^2+2\,n^{-2}\left(\sum_{i=1}^n(Z_i-\mu)/(\sum_{i=1}^n(Z_i-\mu)^2)^{\frac12}\right)^2}{1-n^{-1}\left(\sum_{i=1}^n(Z_i-\mu)/(\sum_{i=1}^n(Z_i-\mu)^2)^{\frac12}\right)^2}$$

$$= \frac{o_P(1)+n^{-2}O_P(1)}{1+n^{-1}O_P(1)} = o_P(1), \quad n\to\infty.$$

Conversely, assume that $\max_{1\le i\le n}(Z_i-\overline{Z})^2/\sum_{i=1}^n(Z_i-\overline{Z})^2 = o_P(1)$, $n\to\infty$. According to Theorem 2, to conclude the FCLT for $\widetilde{T}_n^t(Z_1-\mu,\dots,Z_n-\mu)$, it suffices to show that, as $n\to\infty$, $\max_{1\le i\le n}(Z_i-\mu)^2/\sum_{i=1}^n(Z_i-\mu)^2 = o_P(1)$. As $n\to\infty$, we have

$$\frac{\max_{1\le i\le n}(Z_i-\mu)^2}{\sum_{i=1}^n(Z_i-\mu)^2} \le \frac{2\max_{1\le i\le n}(Z_i-\overline{Z})^2+2(\overline{Z}-\mu)^2}{\sum_{i=1}^n(Z_i-\overline{Z})^2+n(\overline{Z}-\mu)^2}$$

$$\le \frac{2\max_{1\le i\le n}(Z_i-\overline{Z})^2}{\sum_{i=1}^n(Z_i-\overline{Z})^2}+\frac{2}{n} = o_P(1).$$

□

## Acknowledgements

The author is thankful to Miklós Csörgő for his advice and stimulating discussions on the paper.